\documentclass[11pt]{amsart}      
\usepackage{amsfonts}
\usepackage{mathrsfs}
\begin{document}

\title{Left-symmetric Superalgebra Structures on the Super-Virasoro
Algebras}

\author{Xiaoli Kong}
\address{School of Mathematical Sciences, Xiamen
University, Xiamen, Fujian 361005, P.R.
China}\email{kongxl.math@gmail.com}

\author{Chengming Bai}
\address{Chern Institute of
Mathematics and LPMC, Nankai University, Tianjin 300071, P.R.China}
\email{baicm@nankai.edu.cn}

\def\shorttitle{Left-symmetric superalgebra Structures on the Super-Virasoro
Algebras}

\begin{abstract}
In this paper, we classify the compatible left-symmetric
superalgebra structures on the super-Virasoro algebras satisfying
certain natural conditions.

\end{abstract}



\maketitle

\section{Introduction}

Left-symmetric algebras (or under other names like pre-Lie algebras,
quasi-associative algebras, Vinberg algebras and so on) are a class
of natural algebraic systems appearing in many fields in mathematics
and mathematical physics. They were first mentioned by A. Cayley in
1896 as a kind of rooted tree algebras (\cite{Ca}) and arose again
from the study of convex homogenous cones (\cite{V}), affine
manifolds and affine structures on Lie groups (\cite{Ko}),
deformation of associative algebras (\cite{G}) in 1960s. They play
an important role in the study of symplectic and complex structures
on Lie groups and Lie algebras (\cite{AS}, \cite{Ch}, [DaM1-2],
\cite{LM}), phase spaces of Lie algebras (\cite{Ba}, [Ku1-2]),
certain integrable systems (\cite{Bo}, \cite{SS}), classical and
quantum Yang-Baxter equations (\cite{DiM}, \cite{ES}, \cite{GS},
\cite{Ku3}), combinatorics (\cite{E}), quantum field theory
(\cite{CK}), vertex algebras (\cite{BK}), operad (\cite{CL}) and so
on (see a survey article \cite{Bu} and the references therein).

The super-version of left-symmetric algebras, the left-symmetric
superalgebras,  also appeared in a lot of fields ([CL], [G], [VM],
etc.). For example, to our knowledge, they were first introduced by
Gerstenhaber to study the Hochschild cohomology of associative
algebras (\cite{G}).

On the other hand, the Virasoro and super-Virasoro algebras are not
only a class of important infinite-dimensional Lie algebras and Lie
superalgebras, but also one of the principal Lie algebras and Lie
superalgebras of physical interest. For example, they are the
fundamental algebraic structures in conformal and super-conformal
field theory. As it was pointed out in \cite{Ku2}, a compatible
left-symmetric algebra structure on the Virasoro algebra can be
regarded as the ``nature of the Virasoro algebra". In fact, the
compatible left-symmetric algebra on the Virasoro algebra $\mathcal
V$ given in \cite{Ku2} satisfies
$$cc=x_m c=c
x_m=0,\;\;x_mx_n=f(m,n)x_{m+n}+\omega(m,n)c,\eqno (1.1)$$ where
$f(m,n)$ and $\omega(m,n)$ are two complex-value functions, and $\{
x_m, c\mid m\in \mathbb Z\}$ is a basis of the Virasoro algebra
$\mathcal V$ satisfying
$$[c,
x_n]=0,\;\;[x_m,x_n]=(m-n)x_{m+n}+\frac{c}{12}(m^3-m)\delta_{m+n,0}.\eqno
(1.2)$$ The condition (1.1) is natural since it means that the
compatible left-symmetric algebra is still graded and $c$ is also a
central extension given by $\omega(m,n)$.  Moreover, in [KCB], we
proved that any compatible left-symmetric algebra structure on the
Virasoro algebra $\mathcal V$ satisfying equation (1.1) was
isomorphic to one of the examples given in [Ku2].

In this paper, we study the compatible left-symmetric superalgebra
structures on the super-Virasoro algebras. Motivated by the study in
the case of the ordinary Virasoro algebra, we classify such
left-symmetric superalgebras satisfying some natural conditions like
equation (1.1). The paper is organized as follows. In Section 2, we
give some necessary definitions, notations and basic results on
left-symmetric superalgebras and the super-Virasoro algebras. We
also give the classification of compatible left-symmetric algebra
structures on the ordinary Virasoro algebra satisfying equation
(1.1).  In section 3, we study the compatible left-symmetric
superalgebra structures on the centerless super-Virasoro algebras
satisfying certain natural conditions. In section 4, we discuss the
non-trivial central extensions of the left-symmetric superalgebras
obtained in section 3 whose super-commutator is a super-Virasoro
algebra.

Throughout this paper, all algebras are over the complex field
${\mathbb C}$ and the indices $m,n,l\in\mathbb{Z}$ and
$r,s,t\in{\mathbb Z}+\theta$ for $\theta=0$ or $\theta=\frac{1}{2}$,
unless otherwise stated.

\section{Preliminaries and fundamental results}

Let $(A,\cdot)$ be an algebra over a field $\mathbb F$. $A$ is said
to be a superalgebra if the underlying vector space of $A$ is
${\mathbb Z}_2$-graded, i.e., $A=A_{\bar0 }\oplus A_{\bar 1}$, and
$A_\alpha\cdot A_\beta\subset A_{\alpha+\beta}$, for $\alpha,
\beta\in\mathbb{Z}_2$. An element of $A_{\bar0}$ is called even and
an element of $A_{\bar1}$ is called odd.

{\bf Definition 2.1.} A Lie superalgebra is a superalgebra
$A=A_{\bar0 }\oplus A_{\bar 1}$ with an operation $[\ ,\ ]$
satisfying the following conditions:
$$[a,b]=-(-1)^{\alpha\beta}[b,a];\eqno (2.1)$$
$$[a,[b,c]]=[[a,b],c]+(-1)^{\alpha\beta}[b,[a,c]],\eqno (2.2)$$
where $a\in A_\alpha, b\in A_\beta, c\in A, \;\;\alpha, \beta\in
{\mathbb Z}_2$.

{\bf Definition 2.2.} A superalgebra $A$ is called a left-symmetric
superalgebra if the associator
$$(x,y,z):=(x\cdot
y)\cdot z-x\cdot (y\cdot z)\eqno (2.3)$$ of  $A$ satisfies
$$(x,y,z)=(-1)^{\alpha\beta}(y,x,z),\quad\forall\;
x\in A_\alpha,\; y\in A_\beta,\;z\in A,\;\; \alpha, \beta\in{\mathbb
Z}_2.\eqno{(2.4)}$$

Obviously, if $A=A_{\bar0 }\oplus A_{\bar 1}$ is a Lie superalgebra
or a left-symmetric superalgebra, then $A_{\bar 0}$ is an ordinary
Lie algebra or a left-symmetric algebra respectively. Moreover, let
$A$ be a left-symmetric superalgebra, then it is easy to know that
the super-commutator
$$[x,y]=x\cdot y-(-1)^{\alpha\beta}y\cdot x,\quad \forall\; x\in A_\alpha,\
y\in A_\beta,\ \alpha,\beta\in{\mathbb Z}_2, \eqno(2.5)$$ defines a
Lie superalgebra ${\mathcal G}(A)$ which is called the sub-adjacent
Lie superalgebra of $A$ and $A$ is also called the compatible
left-symmetric superalgebra structure on the Lie superalgebra
${\mathcal G}(A)$.

On the other hand, we recall the definition of the super-Virasoro
algebras. There are two super-Virasoro algebras which correspond to
$N=1$ ([R]) and $N=2$ ([NS1-2]) super-conformal field theory
respectively. In fact, let $\theta=0$ or $\displaystyle\frac{1}{2}$
which corresponds to the Ramond case ([R]) or the Neveu-Schwarz case
([NS1-2]) respectively. Let ${\mathcal S \mathcal V}={\mathcal
S\mathcal V}_{\bar0}\oplus {\mathcal S \mathcal V}_{\bar1}$ denote a
super-Virasoro algebra with a basis $\{L_m,G_r, c\mid m\in{\mathbb
Z},r\in{\mathbb Z}+\theta\}$. The super-brackets are defined as
follows,
$$\left.\begin{split}
&[L_m,L_n]=(m-n)L_{m+n}+\displaystyle\frac{c}{12}(m^3-m)\delta_{m+n,0},\\
&[L_m,G_r]=\displaystyle(\frac{m}{2}-r)G_{m+r},\\
&[G_r,G_s]=2L_{r+s}+\displaystyle\frac{c}{12}(4r^2-1)\delta_{r+s,0},\\
&[{\mathcal S \mathcal V}_{\bar0},c]=[{\mathcal S \mathcal
V}_{\bar1},c]=0,
\end{split}\right.\eqno{(2.6)}$$
where the even subspace ${\mathcal S \mathcal V}_{\bar0}$ is spanned
by $\{L_m, c\mid m\in \mathbb Z\}$ and the odd subspace ${\mathcal
S\mathcal V}_{\bar1}$ is spanned by $\{G_r\mid r\in \mathbb
Z+\theta\}$. Obviously, ${\mathcal S\mathcal V}_{\bar0}$ is nothing
but an ordinary Virasoro algebra. A class of compatible
left-symmetric algebra structures on the ordinary Virasoro algebra
satisfying equation (1.1) were given in [Ku2]. Moreover, such
left-symmetric algebras were classified in [KCB].

{\bf Theorem 2.3.} [KCB] Any compatible left-symmetric algebra
structure on the Virasoro algebra ${\mathcal S \mathcal V}_{\bar0}$
satisfying equation (1.1) is isomorphic to one of the (mutually
non-isomorphic) left-symmetric algebras given by the multiplication
$$L_mL_n=\frac{-n(1+\epsilon
n)}{1+\epsilon(m+n)}L_{m+n}+\frac{c}{24}(m^3-m+(\epsilon-\epsilon^{-1})m^2)\delta_{m+n,0},\
\forall\; m,n\in\mathbb{Z},\eqno(2.7)$$ where $m,n\in\mathbb{Z}$,
$c$ is an annihilator and $\rm{Re} \epsilon>0,
\epsilon^{-1}\notin\mathbb{Z}\ \mbox{ or }\ \rm{Re}
\epsilon=0,\rm{Im} \epsilon>0.$

\section{Compatible left-symmetric superalgebra structures on the
centerless super-Virasoro algebras}

Let $\widetilde{\mathcal S \mathcal V}=\widetilde{\mathcal S\mathcal
V}_{\bar0}\oplus \widetilde{\mathcal S \mathcal V}_{\bar1}$ be a
centerless super-Virasoro algebra with a basis $\{L_m,G_r\mid
m\in{\mathbb Z},r\in{\mathbb Z}+\theta\}$ and the super-brackets be
given in equation (2.6) with $c=0$. Motivated by Theorem 2.3, it is
natural to consider the compatible left-symmetric superalgebra
structures on $\widetilde{\mathcal S \mathcal V}$ also satisfy the
``graded" condition, that is, the multiplications of the compatible
left-symmetric superalgebra structures on $\widetilde{\mathcal S
\mathcal V}$ satisfy
$$\left.\begin{split}
&L_m\cdot L_n=f(m,n)L_{m+n},\quad L_m\cdot
G_r=g(m,r)G_{m+r},\\
&G_r\cdot L_m=h(r,m)G_{m+r},\quad G_r\cdot
G_s=d(r,s)L_{r+s},\end{split}\right.\eqno(3.1)$$ where $f$, $g$, $h$
and $d$ are $\mathbb{C}$-value functions. Then the super-commutators
give the  super-Virasoro algebra $\widetilde{\mathcal S \mathcal V}$
if and only if $f(m,n),g(m,r),h(r,m),d(r,s)$ satisfy
$$f(m,n)-f(n,m)=m-n,\quad g(m,r)-h(r,m)=\frac{m}{2}-r,\quad
d(r,s)+d(s,r)=2.\eqno(3.2)$$ Furthermore, the functions $f(m,n),\
g(m,r),\ h(r,m)$ and $d(r,s)$ can define a left-symmetric
superalgebra with a basis $\{L_m,G_r\mid m\in{\mathbb
Z},r\in{\mathbb Z}+\theta\}$ if and only if they satisfy the
following equations,
$$(L_m,L_n,L_l)=(-1)^{0\cdot0}(L_n,L_m,L_l),\quad
(L_m,L_n,G_r)=(-1)^{0\cdot0}(L_n,L_m,G_r),$$
$$(L_m,G_r,L_n)=(-1)^{0\cdot1}(G_r,L_m,L_n),\quad
(L_m,G_r,G_s)=(-1)^{0\cdot1}(G_r,L_m,G_s),$$
$$(G_r,G_s,L_m)=(-1)^{1\cdot1}(G_s,G_r,L_m),\quad
(G_r,G_s,G_t)=(-1)^{1\cdot1}(G_s,G_r,G_t).$$ The above equations are
equivalent to the following equations,
$$\left\{\begin{split}
&(m-n)f(m+n,l)=f(n,l)f(m,n+l)-f(m,l)f(n,m+l),\\
&(m-n)g(m+n,r)=g(n,r)g(m,n+r)-g(m,r)g(n,m+r),\\
&(\frac{m}{2}-r)h(m+r,n)=h(r,n)g(m,n+r)-f(m,n)h(r,m+n),\\
&(\frac{m}{2}-r)d(m+r,s)=d(r,s)f(m,r+s)-g(m,s)d(r,m+s),\\
&2f(r+s,m)=h(s,m)d(r,m+s)+h(r,m)d(s,m+r),\\
&2g(r+s,t)=d(s,t)h(r,s+t)+d(r,t)h(s,r+t).
\end{split}\right.\eqno(3.3)$$

{\bf Proposition 3.1.} Any compatible left-symmetric superalgebra
structure $\widetilde{V}$ on $\widetilde{\mathcal S\mathcal V}$
satisfies equation (3.1) if and only if the functions in equation
(3.1) satisfy equations (3.2) and (3.3).

By Theorem 2.3, we only need to consider the case that
$$f(m,n)=\displaystyle\frac{-n(1+\epsilon n)}{1+\epsilon(m+n)},\eqno (3.4)$$
where $\rm{Re} \epsilon>0, \epsilon^{-1}\notin\mathbb{Z}\ \mbox{ or
}\ \rm{Re} \epsilon=0,\rm{Im} \epsilon>0.$

{\bf Theorem 3.2.} For a fixed $\epsilon$ satisfying $\rm{Re}
\epsilon>0, \epsilon^{-1}\notin\mathbb{Z}\ \mbox{ or }\ \rm{Re}
\epsilon=0,\rm{Im} \epsilon>0$ and $f(m,n)$ satisfying equation
(3.4), there is exactly one solution satisfying equations (3.2) and
(3.3) given by
$$g(m,r)=\displaystyle\frac{-(\frac{m}{2}+r)(1+2\epsilon
r)}{1+2\epsilon(m+r)},\ h(r,m)=\displaystyle\frac{-m(1+\epsilon
m)}{1+2\epsilon(m+r)},\ d(r,s)=\displaystyle\frac{1+2\epsilon
s}{1+\epsilon(r+s)},\eqno (3.5)$$ for
$m,n\in\mathbb{Z},r,s\in\mathbb{Z}+\theta$, which define a
compatible left-symmetric superalgebra ${\widetilde V}_\epsilon$ on
$\widetilde{\mathcal S\mathcal V}$.

{\bf Proof.} It is easy to verify that $f(m,n)$ given in equation
(3.4) and $g(m,r), h(r,m), d(r,s)$ given in equation (3.5) satisfy
equations (3.2)--(3.3). On the other hand, set
$$ G(m,r)=g(m,r)\frac{1+2\epsilon(m+r)}{1+2\epsilon r},$$$$
H(r,m)=h(r,m)\frac{1+2\epsilon(m+r)}{1+\epsilon m},$$$$
D(r,s)=d(r,s)\frac{1+\epsilon(r+s)}{1+2\epsilon s}.$$ Then we only
need to prove that $$G(m,r)=-\frac{m}{2}-r,\quad H(r,m)=-m,\quad
D(r,s)=1.$$ We rewrite equations (3.2) and (3.3) involving $g(m,r),
h(r,m), d(r,s)$ as
$$G(m,r)(1+2\epsilon r)-H(r,m)(1+\epsilon
m)=(\frac{m}{2}-r)(1+2\epsilon(m+r)),\eqno(3.6)$$
$$D(r,s)(1+2\epsilon s)+D(s,r)(1+2\epsilon r)=2+2\epsilon (r+s),\eqno(3.7)$$
$$(m-n)G(m+n,r)=G(n,r)G(m,n+r)-G(m,r)G(n,m+r),\eqno(3.8)$$
$$(\frac{m}{2}-r)H(m+r,n)=H(r,n)G(m,n+r)+nH(r,m+n),\eqno(3.9)$$
$$(\frac{m}{2}-r)D(m+r,s)=-(r+s)D(r,s)-G(m,s)D(r,m+s),\eqno(3.10)$$
$$-2m=H(s,m)D(r,m+s)+H(r,m)D(s,m+r),\eqno(3.11)$$
$$2G(r+s,t)=D(s,t)H(r,s+t)+D(r,t)H(s,r+t). \eqno(3.12)$$
Let $r=s$ in equation $(3.7)$, we have $$D(s,s)=1,\quad \forall\;
s\in\mathbb{Z}+\theta.\eqno(3.13)$$ In fact, $D(r,s)\neq0$ for all
$r,s\in\mathbb{Z}+\theta$. Otherwise, assume there exist $r_1,s_1$,
such that $D(r_1,s_1)=0$. Let $r=s=r_1, m=s_1-r_1$ in equation
$(3.11)$, we have
$$-(s_1-r_1)=H(r_1,s_1-r_1)D(r_1,s_1)=0.$$
Hence $r_1=s_1$. It is contradictory  to equation $(3.13)$.

Let $m=0, r=s$ in equations $(3.10)$ and $(3.6)$, $r=s=-t$ in
equation $(3.12)$, $m=-2s,r=s$ in equations $(3.6)$ and $(3.10)$,
and $m=-2s\neq0, r=3s$ in  equation $(3.10)$ respectively, we know
that
$$\left.\begin{split} &G(0,s)=-s,\ H(s,0)=0,\ G(2s,-s)=0,\
H(s,-2s)=2s,&\\
&D(-s,s)=D(3s,s)=1,\quad \forall\;
s\in\mathbb{Z}+\theta.&\end{split}\right.\eqno(3.14)$$ Let
$m=-2(n+r)$ in equation $(3.9)$ and $m+r+s=0$ in equation $(3.11)$,
we have
$$\left\{\begin{split}
&-(n+2r)H(-2n-r,n)=nH(r,-n-2r),\\
&-2m=H(-m-r,m)+H(r,m).\end{split}\right.\eqno(3.15)$$ Let $r=s$ in
equation $(3.11)$, we have
$$-m=H(r,m)D(r,m+r).$$ So
$$H(r,m)=\displaystyle\frac{-m}{D(r,m+r)}.\eqno(3.16)$$
By equations $(3.13)$, $(3.14)$, $(3.15)$, and $(3.16)$, we have
$$\left\{\begin{split}
&D(-2n-r,-n-r)=D(r,-n-r),\\
&\frac{1}{D(r,m+r)}+\frac{1}{D(-m-r,-r)}=2.\end{split}\right.\eqno(3.17)$$
Let $-n-r=s, m+r=s$ in equations $(3.17)$, we have
$$D(r,s)=D(2s+r,s),\;\;{\rm  and}\;\;D(-s,-r)=D(-s,-2s-r),\quad\forall\; r,s\in\mathbb{Z}+\theta.$$
Thus by induction, we know that
$$D(r,s)=D(2ks+r,s),\;\;D(-s,-r)=D(-s,-2ks-r),\;\;\forall\;k\in\mathbb{Z}.\eqno(3.18)$$
Therefore, we have
$$D(r,r)=D((2k+1)r,r)=D(r,(2k+1)r)=1,\quad \forall\; k\in\mathbb{Z}.$$
Let $r=s$ in equation $(3.12)$. Then by equations $(3.16)$ and
$(3.18)$, we have
$$G(2s,t)=D(s,t)H(s,s+t)=D(s,t)\frac{-s-t}{D(s,2s+t)}=-s-t.$$
Let $m=2s$ in equation $(3.6)$, we have
$$H(t,2s)=-2s,\quad\forall\;s,t\in\mathbb{Z}+\theta.\eqno(3.19)$$

{\it Case (I)}\ $\theta=\displaystyle\frac{1}{2}$. Then
$$D(\theta,\pm\theta)=D(\pm\theta,\theta)=1.$$ Hence
$$D(k+\theta,\pm\theta)=D(\pm\theta,k+\theta)=1.$$ That is,
$$\displaystyle D(r,\pm\frac{1}{2})=D(\pm\frac{1}{2},r)=1,
\;\;\forall\; r\in\mathbb{Z}+\theta.$$ Assume that for any
$|r_1|\leq|s_1|$, we have $D(r_1,s_1)=1$. Then
$$D(r_1,s_1)=D(2ks_1+r_1,s_1)=1,\;\;{\rm
and}\;\;D(s_1,r_1)=D(s_1,2ks_1+r_1)=1.$$ For any
$r\in\mathbb{Z}+\theta,$ there exist
$k\in\mathbb{Z},\;r_1\in\mathbb{Z}+\theta,$ and $|r_1|\leq|s_1|$,
such that $r=2ks_1+r_1$. Therefore,
$$D(r,s_1)=D(s_1,r)=1,\;\;\forall\; r\in\mathbb{Z}+\theta.$$ Hence by
induction, we know that $D(r,s)=1$ for any
$r,s\in\mathbb{Z}+\theta$. Therefore,
$$H(r,m)=-m,\quad G(m,r)=-\frac{m}{2}-r,\quad \forall\; m\in\mathbb{Z},
r\in\mathbb{Z}+\theta.$$

{\it Case (II)}\ $\theta=0$. Let $m=-2t\neq0, s=r=2t$ in equation
$(3.11)$. Then by equation $(3.19)$, we have
$$2t=H(2t,-2t)D(2t,0)=2tD(2t,0).$$
Therefore, we have $$D(2t,0)=1,\quad
D(0,2t)=1,\;\;\forall\;t\in\mathbb{Z}.\eqno(3.20)$$ Let
$r=0,m=s\neq0$ in equations $(3.10)$ and $(3.11)$, we have
$$\left\{\begin{split}
&\frac{m}{2}=-mD(0,m)-G(m,m),\\
&-2m=H(m,m)+H(0,m).\end{split}\right.\eqno(3.21)$$ So
$$H(m,m)=-2m-H(0,m)=-2m+\frac{m}{D(0,m)}=-2m-\frac{2m^2}{m+2G(m,m)}.$$
By equation $(3.6)$, we have $$H(m,m)=-m, \
G(m,m)=-\frac{3m}{2},$$ or
$$H(m,m)=\displaystyle\frac{-m}{1+\epsilon m},\
G(m,m)=\displaystyle\frac{-3m-4\epsilon m^2}{2+4\epsilon m}.$$ In
fact, there does not exist the latter case for any $m\neq0$.
Otherwise, assume that there exists a nonzero integer $m_1$, such
that
$$H(m_1,m_1)=\frac{-m_1}{1+\epsilon m_1},\
G(m_1,m_1)=\frac{-3m_1-4\epsilon m_1^2}{2+4\epsilon m_1}.$$ Then
$$D(0,m_1)=-\frac{1}{2}-\frac{G(m_1,m_1)}{m_1}
=\frac{1+\epsilon m_1}{1+2\epsilon m_1}.$$ Let $m=-s=m_1$, $r=0$ or
$-m_1$ in equation $(3.10)$, we have
\begin{eqnarray*}\left\{\begin{split}
&\frac{m_1}{2}=m_1D(0,-m_1)-G(m_1,-m_1),\\
&\frac{3m_1}{2}D(0,-m_1)=2m_1-G(m_1,-m_1)D(-m_1,0).\end{split}\right.\end{eqnarray*}
Hence
$$\displaystyle\frac{3}{2}D(0,-m_1)=2-D(0,-m_1)D(-m_1,0)+\frac{1}{2}D(-m_1,0).$$
By equation $(3.7)$, we have $$D(0,-m_1)=1,\ D(-m_1,0)=1,$$ or
$$D(0,-m_1)=\frac{3-\epsilon m_1}{1-2\epsilon
m_1},\ D(-m_1,0)=-1-\epsilon m_1.$$ Since $\epsilon\neq0,
\epsilon^{-1}\notin\mathbb{Z}$, we know that
$$\frac{1}{D(0,m_1)}+\frac{1}{D(-m_1,0)}\neq2,$$
which is contradictory to equation $(3.17)$. Hence
$$H(m,m)=-m, \;\;G(m,m)=-\displaystyle\frac{3m}{2}, \quad\forall\;
m\in\mathbb{Z}.$$ By equations $(3.21)$ and $(3.6)$, we have
$$H(0,m)=-m,\
G(m,0)=-\displaystyle\frac{m}{2}.$$ Let $r=0,$ and $m\neq0$ in
equations $(3.8)$ and $(3.9)$, we have
\begin{eqnarray*}\left\{\begin{split}
&n^2-m^2=-nG(m,n)+mG(n,m),\\
&\frac{m}{2}H(m,n)=-nG(m,n)-n(m+n).\end{split}\right.\end{eqnarray*}
So
$$H(m,n)+2G(n,m)=-2(m+n).$$ By equations $(3.6), (3.13)$ and $(3.16)$,  we know that
$$H(m,n)=-n, \quad G(n,m)=-\frac{n}{2}-m, \quad D(m,n)=1, \quad
\forall\; m,n\in\mathbb{Z}.\eqno \Box$$

\section{Compatible left-symmetric superalgebra structures on the
super-Virasoro algebras}

In this section, we consider the central extensions of the
left-symmetric superalgebras obtained in Section 3 whose
super-commutator is a super-Virasoro algebra $\mathcal S \mathcal
V$.

Let $\tilde{A}$ be a left-symmetric superalgebra and
$\omega:\tilde{A}\times \tilde{A}\rightarrow \mathbb C$ be a
bilinear form. It defines a multiplication on the space
$A=\tilde{A}\oplus {\mathbb C} c$, by the rule
$$(x+\lambda c)\cdot (y+\mu c)=x\cdot y +\omega(x,y)c,\quad \forall\; x,y\in A,\ \lambda,\mu\in \mathbb C.\eqno{(4.1)}$$
Let $$B(x,y,z):=\omega(x\cdot y, z)-\omega(x, y\cdot z).\eqno
(4.2)$$ Then it is easy to know that $A$ is a left-symmetric
superalgebra if and only if
$$B(x,y,z)=(-1)^{\alpha\beta}B(y,x,z),\quad\forall
\;x\in \tilde{A}_\alpha,\ y\in \tilde{A}_\beta,\ z\in
\tilde{A},\quad\alpha,\beta\in{\mathbb Z}_2,\eqno{(4.3)}$$ $ A$ is
called a central extension of $\tilde A$. Moreover, by construction,
the bilinear form
$$\Omega(x,y)=\omega(x,y)-(-1)^{\alpha\beta}\omega(y,x),\quad\forall\;
 x\in \tilde{A}_\alpha,\ y\in \tilde{A}_\beta,z\in\tilde{A},\quad\alpha,\beta\in{\mathbb Z}_2,\eqno{(4.4)}$$
defines a central extension of the Lie superalgebra ${\mathcal
G}(A)$.

Let the left-symmetric superalgebra $\widetilde{V}_\epsilon$ on a
centerless Virasoro algebra $\widetilde {\mathcal S\mathcal V}$ be
given through Theorem 3.2. Since a super-Virasoro algebra $\mathcal
S \mathcal V$ is a central extension of a centerless super-Virasoro
algebra $\widetilde {\mathcal S\mathcal V}$, it is natural to
consider the central extension $V_\epsilon=\widetilde
{V}_\epsilon\oplus \mathbb C c$ of $\widetilde{V}_\epsilon$ such
that $V_\epsilon$ is a compatible left-symmetric superalgebra
structure on the super-Virasoro algebra $\mathcal S\mathcal V$ with
$c$ being the annihilator of $V_\epsilon$, that is, the products of
$V_\epsilon$ are given by $$\left.\begin{split}
&L_m\cdot L_n=f(m,n)L_{m+n}+\omega(L_m,L_n)c,\\
&L_m\cdot G_r=g(m,r)G_{m+r}+\omega (L_m, G_r)c,\\
&G_r\cdot L_m=h(r,m)G_{m+r}+\omega(G_r, L_m)c,\\
&G_r\cdot G_s=d(r,s)L_{r+s}+\omega(G_r,G_s)c,\\
&c\cdot c=c\cdot L_m=L_m\cdot c=c\cdot G_r=G_r\cdot c=0,
\end{split}\right. \eqno(4.5)$$
where the functions $f(m,n), g(m,r), h(r,m)$ and $d(r,s)$ satisfy
equations (3.4) and (3.5).

For convenience, set
$$\left.\begin{split}&\omega(L_m,L_n)=\varphi(m,n),\quad \omega(L_m,G_r)=\psi(m,r),\\
&\omega(G_r, L_m)=\rho(r,m),\quad
\omega(G_r,G_s)=\sigma(r,s).\end{split}\right.\eqno{(4.6)}$$ So the
super-commutators of $V_\epsilon$ give a super-Virasoro algebra
${\mathcal S \mathcal V}$ if and only if $\varphi(m,n)$,
$\psi(m,r)$, $\rho(r,m)$, and $\sigma(r,s)$ satisfy
$$\left.\begin{split}
&\varphi(m,n)-\varphi(n,m)=\frac{1}{12}(m^3-m)\delta_{m+n,0},\\
&\sigma(r,s)+\sigma(s,r)=\frac{1}{12}(4r^2-1)\delta_{r+s,0},\\
&\psi(m,r)-\rho(r,m)=0.\end{split}\right.\eqno(4.7)$$ By equation
(4.3), we have
$$B(L_m,L_n,L_l)=(-1)^{0\cdot0}B(L_n,L_m,L_l),\;\
B(L_m,L_n,G_r)=(-1)^{0\cdot0}B(L_n,L_m,G_r),$$
$$B(L_m,G_r,L_n)=(-1)^{0\cdot1}B(G_r,L_m,L_n),\;\
B(L_m,G_r,G_s)=(-1)^{0\cdot1}B(G_r,L_m,G_s),$$
$$B(G_r,G_s,L_m)=(-1)^{1\cdot1}B(G_s,G_r,L_m),\;\
B(G_r,G_s,G_t)=(-1)^{1\cdot1}B(G_s,G_r,G_t).$$ They are equivalent
to the following equations
$$(m-n)\varphi(m+n,l)=f(n,l)\varphi(m,n+l)-f(m,l)\varphi(n,m+l),\eqno(4.8)$$
$$(m-n)\psi(m+n,r)=g(n,r)\psi(m,n+r)-g(m,r)\psi(n,m+r),\eqno(4.9)$$
$$(\frac{m}{2}-r)\rho(m+r,n)=h(r,n)\psi(m,n+r)-f(m,n)\rho(r,m+n),\eqno(4.10)$$
$$(\frac{m}{2}-r)\sigma(m+r,s)=d(r,s)\varphi(m,r+s)-g(m,s)\sigma(r,m+s),\eqno(4.11)$$
$$2\varphi(r+s,m)=h(s,m)\sigma(r,m+s)+h(r,m)\sigma(s,m+r),\eqno(4.12)$$
$$2\psi(r+s,t)=d(s,t)\rho(r,s+t)+d(r,t)\rho(s,r+t).\eqno(4.13)$$

{\bf Proposition 4.1.} Any compatible left-symmetric superalgebra
structure $V$ on $\mathcal S\mathcal V$ satisfies equation (4.5) if
and only if the functions in equation (4.5) satisfy equations (3.4),
(3.5) and (4.7)--(4.13).

If a central extension $V_{\epsilon}$ of $\widetilde{V}_{\epsilon}$
given by $\omega$ satisfying equation (4.5) defines a compatible
left-symmetric superalgebra structure on $\mathcal S\mathcal V$,
then $\varphi(m,n)$ defines a central extension of
$\widetilde{\mathcal{SV}}_{\bar0}$. By Theorem 2.3, we know that
$$\varphi(m,n)=\displaystyle\frac{1}{24}(m^3-m+(\epsilon-\epsilon^{-1})m^2)\delta_{m+n,0}.\eqno
(4.14)$$

{\bf Theorem 4.2.} For a fixed $\epsilon \in \mathbb C$ satisfy
$\rm{Re} \epsilon>0, \epsilon^{-1}\notin\mathbb{Z}\ \mbox{ or }\
\rm{Re} \epsilon=0,\rm{Im} \epsilon>0$, let the functions $f(m,n),
g(m,r), h(r,m)$ and $d(r,s)$ satisfy equations (3.4) and (3.5), and
$\varphi(m,n)$ satisfy equation (4.14). Then there is exactly one
solution satisfying equations (4.7)--(4.13) given by
$$\sigma(r,s)=\frac{1}{24}(4r^2-1+2(\epsilon-\epsilon^{-1})r)\delta_{r+s,0},\
\phi(m,r)=\rho(r,m)=0,\eqno (4.15)$$ for
$m\in\mathbb{Z},r,s\in\mathbb{Z}+\theta,$ which define a compatible
left-symmetric superalgebra ${V}_\epsilon$ on ${\mathcal S\mathcal
V}$.

{\bf Proof.} It is easy to verify that $\varphi(m,n)$ given in
equation (4.14) and $\sigma(r,s),\phi(m,r),\rho(r,m)$ given in
equation (4.15) satisfy equations (4.7 )--(4.13).

On the other hand, let $m=0, r+s\neq0$ in equation (4.11). Then we
have
$$-r\sigma(r,s)=d(r,s)\varphi(0,r+s)-g(0,s)\sigma(r,s)=s\sigma(r,s).$$
Hence
$$\sigma(r,s)=0,\;\;\forall\; r+s\neq0.$$
Let $r=s,m=-2s$ in equation (4.12), we have
$$2\varphi(2s,-2s)=h(s,-2s)\sigma(s,-s)+h(s,-2s)\sigma(s,-s)=4s\sigma(s,-s).$$
So
$$\displaystyle\sigma(s,-s)=\frac{1}{24}(4s^2-1+2(\epsilon-\epsilon^{-1})s).$$
Thus
$$\sigma(r,s)=\frac{1}{24}(4r^2-1+2(\epsilon-\epsilon^{-1})r)\delta_{r+s,0},\quad
\forall\; r,s\in\mathbb{Z}+\theta.$$

Next, we prove that
$$\psi(m,r)=\rho(r,m)=0,\;\;
\forall\; m\in\mathbb{Z},\ r\in\mathbb{Z}+\theta.$$ There are two
cases as follows.

{\it  Case (I)}\quad $\theta=\displaystyle\frac{1}{2}$. Let $m=n=0$
in equation (4.10), we have
$$-r\rho(r,0)=h(r,0)\psi(0,r)-f(0,0)\rho(r,0)=0.$$ So $\rho(r,0)=0$.
By equation (4.7), we know that $\psi(0,r)=0$. Let $n=0$ in equation
(4.9), we have
$$m\psi(m,r)=\psi(m,r)g(0,r)-\psi(0,m+r)g(m,r)=-r\psi(m,r).$$ Hence
$(m+r)\psi(m,r)=0$. Therefore, we have
$$\psi(m,r)=\rho(r,m)=0,\;\;\forall\; m\in\mathbb{Z},
r\in\mathbb{Z}+\frac{1}{2}.$$

{\it Case (II)}\quad $\theta=0$. Let $n=0,m=-r\neq0$ in equation
$(4.10)$, we have $$\psi(0,0)=\rho(0,0)=0.$$ Let $m=n=0, r\neq0$ in
equation (4.10) and $m=r=0, n\neq0$ in equation $(4.9)$
respectively, we have
$$\rho(r,0)=\psi(0,r)=0,\;\;\psi(n,0)=\rho(0,n)=0,\;\;\forall\; r,n\in\mathbb{Z}, r,n\neq0.$$
Let $r=0$, $m,n\neq0$ in equations (4.9) and (4.10), we have
\begin{eqnarray*}\left\{\begin{split}
&\psi(m,n)\frac{n}{1+2\epsilon n}-\psi(n,m)\frac{m}{1+2\epsilon m}=0,\\
&\frac{m}{2}\psi(n,m)+\psi(m,n)\frac{n(1+\epsilon n)}{1+2\epsilon
n}=0.\end{split}\right.\end{eqnarray*} Since
$\epsilon\neq0,\epsilon^{-1}\notin\mathbb{Z}$, we have
$$\psi(n,m)=\rho(n,m)=0,\;\; \forall\; m,n\in\mathbb{Z}.\eqno\Box$$

By Theorem 2.3, it is easy to know that $V_\epsilon$ are mutually
non-isomorphic for all $\epsilon\in \mathbb C$ satisfying $\rm{Re}
\epsilon>0, \epsilon^{-1}\notin\mathbb{Z}\ \mbox{ or }\ \rm{Re}
\epsilon=0,\rm{Im} \epsilon>0$. Furthermore, by Theorem 2.3,
Proposition 3.1, Theorem 3.2, Proposition 4.1 and Theorem 4.2
together, we have the following conclusion.

{\bf Theorem 4.3.} Any compatible left-symmetric superalgebra on a
super-Virasoro algebra satisfying equation (4.5) is isomorphic to
one of the following (mutually non-isomorphic) left-symmetric
superalgebras given by the multiplication
\begin{eqnarray*}
L_m\cdot L_n&=&\frac{-n(1+\epsilon
n)}{1+\epsilon(m+n)}L_{m+n}+\frac{c}{24}(m^3-m+(\epsilon-\epsilon^{-1})
m^2)\delta_{m+n,0},\\
L_m\cdot G_r&=&\frac{-(\frac{m}{2}+r)(1+2\epsilon r)}{1+2\epsilon(m+r)}G_{m+r},\\
G_r\cdot L_m&=&\frac{-m(1+\epsilon m)}{1+2\epsilon(m+r)}G_{m+r},\\
G_r\cdot G_s&=&\frac{1+2\epsilon
s}{1+\epsilon(r+s)}L_{r+s}+\frac{c}{24}(4r^2-1+2(\epsilon-\epsilon^{-1})r)\delta_{r+s,0},
\end{eqnarray*}
where $m,n\in\mathbb{Z},r,s\in\mathbb{Z}+\theta$, $c$ is  an
annihilator and $\rm{Re} \epsilon>0, \epsilon^{-1}\notin\mathbb{Z}\
\mbox{ or }\ \rm{Re} \epsilon=0,\rm{Im} \epsilon>0.$

\section*{Acknowledgements}  The second author thanks Professor B.A. Kupershmidt
for important suggestion and encouragement. This work was supported
by the National Natural Science Foundation of China (10571091,
10621101), NKBRPC (2006 CB805905), Program for New Century Excellent
Talents in University.

\end{document}